\RequirePackage{fix-cm}

\documentclass{article}

\usepackage{amsmath}
\usepackage{amsfonts}

\usepackage{amssymb}
\usepackage{latexsym}
\usepackage{bbold}		
\usepackage{url}		
\usepackage{hyperref}		

\usepackage{mathtools}		

\usepackage[font=small,labelfont=bf]{caption}

\usepackage{mathptmx}       
\usepackage{helvet}         
\usepackage{courier}        
\usepackage{type1cm}        
\usepackage{makeidx}         
\usepackage{graphicx}        
\usepackage{multicol}        
\usepackage[bottom]{footmisc}
\usepackage{xspace}
\usepackage{bbm}

\usepackage{amsmath,amssymb}

\usepackage{draftwatermark}
\SetWatermarkText{Preprint}
\SetWatermarkScale{1}

\newcommand{\phispace}{V_h}
\newcommand{\comment}[1]{}

\newcommand{\grid}{{\mathcal{T}_h}}
\newcommand{\entity}{K}
\newcommand{\elem}{\entity}
\newcommand{\isec}{e}

\newcommand{\NNN}{{\mathbbm N}}

\newcommand{\maxnuminterface}{\mathcal{N}_\grid}
\newcommand{\nspec}{n_{spec}}

\newcommand{\flux}[1]{\widehat{#1}}

\newcommand{\fluxF}{\flux{\mathcal{F}}_{\isec}}

\newcommand{\fluxA}{\flux{\mathcal{A}}_{\isec}}

\newcommand{\vect}[1]{\boldsymbol{#1}}
\newcommand{\vecU}{\vect{U}}
\newcommand{\vecV}{\vect{V}}

\newcommand{\vecv}{\vect{v}}
\newcommand{\vecw}{\vect{w}}
\newcommand{\basefct}{\vect{\varphi}}

\newcommand{\nbold}{\boldsymbol{n}}
\newcommand{\taubold}{\boldsymbol{\tau}}
\newcommand{\df}{\vecU_h}

\newcommand{\oper}[1]{\mathcal{#1}}
\newcommand{\spcoper}{\oper{L}_h}

\newcommand{\vjump}[1]{[ \! [ {#1} ] \! ]_{\isec} }

\newcommand{\vaver}[1]{\{ \! \! \{ {#1} \} \! \! \}_{\isec} }

\newcommand{\su}{S(\df)}

\newcommand{\liftr}{{\boldsymbol{r}}}
\newcommand{\liftre}{{\liftr_e}}
\newcommand{\liftfactor}{\chi}
\newcommand{\Kminus}{{K^-_e}}
\newcommand{\dune}{\textsc{Dune}\xspace}
\newcommand{\dunefem}{\textsc{Dune-Fem}\xspace}

\makeindex             

\begin{document}

\title{Efficient Parallel Simulation of Atherosclerotic Plaque Formation Using Higher Order Discontinuous Galerkin Schemes}


\author{Stefan Girke\thanks{Contact: \href{mailto:stefan.girke@wwu.de}{\nolinkurl{stefan.girke@wwu.de}}, \href{mario.ohlberger@wwu.de}{\nolinkurl{mario.ohlberger@wwu.de}}, Institute for Computational and Applied Mathematics at the University of M{\"u}nster, Einsteinstra{\ss}e 62, D-48149 M{\"u}nster, Germany} \and Robert Kl\"{o}fkorn\thanks{Contact: \href{robertk@ucar.edu}{\nolinkurl{robertk@ucar.edu}}, Center for Atmospheric Research, 1850 Table Mesa Drive, Boulder, CO 80305, USA} \and Mario Ohlberger\footnotemark[1]{}}

\date{}

\maketitle

\begin{abstract}
The compact Discontinuous Galerkin 2 (CDG2) method was successfully tested for elliptic problems, scalar convection-diffusion equations and compressible Navier-Stokes equations. In this paper we use the newly developed DG method to solve a mathematical model for early stages of atherosclerotic plaque formation.
Atherosclerotic plaque is mainly formed by accumulation of lipid-laden cells in the arterial walls which leads to a heart attack in case the artery is occluded or a thrombus is built through a rupture of the plaque.
After describing a mathematical model and the discretization scheme, we present some benchmark tests comparing the CDG2 method to other commonly used DG methods. Furthermore, we take parallelization and higher order discretization schemes into account.
\end{abstract}

\section{Introduction}

Atherosclerotic plaque formation is today seen as a chronic inflammation of the arterial wall which grows over decades and may finally lead to a heart attack in case the artery is occluded or a thrombus is built through a rupture of the plaque. To understand the mechanisms of the chronic inflammation it was recently shown in \cite{kuhlmann2012implantation} that genetically modified (apoE knockdown) mice with a cuff around their carotid develop atherosclerotic plaque formation up- and downstream of the cuff after they were fed with a Western diet. A low wall shear stress of the blood onto the arterial wall or highly oscillating blood flow was shown to be an important indicator for the development of plaque because it damages the endothelial layer. 

At this point our mathematical model (cf. \cite{ibragimov2005mathematical}) comes into play which we want to present in section 2: A dysfunction of the endothelial allows low-density lipoproteins (LDL) to enter the artery wall. Once inside the arterial wall, the LDL becomes oxidized which leads to a recruitment of immune cells, i.e. monocytes. Monocytes differentiate into active macrophages when inside the arterial wall starting continuously absorbing the oxidized LDL. Finally, the macrophages differentiate into foam cells, die and build a necrotic core. Smooth muscle cells (SMCs) from the outer regions of the arterial wall can migrate into the lesion and either become an apoptotic cell or migrate around the lesion to form a fibromuscular cap overlaying the plaque.   

Section 3 describes the spatial and temporal discretization of the CDG2 method which was successfully tested for elliptic problems, scalar convection-diffusion equations and compressible Navier-Stokes equations in \cite{brdar2012cdg2,klofkorn2011benchmark,dgimpl:12}.

We summarize our paper with some 2D and 3D benchmark tests. in section 4 and a conclusion in section 5.

\section{Mathematical Model for Atherosclerotic Inflammation}

A variety of mathematical models dealing with atherosclerotic plaque formation exist, 
see  \cite{calvez2010mathematical,ibragimov2005mathematical}. Here, we focus on six species: immune cells $n_1$ (we do not distinguish between monocytes and macrophages, here), SMCs $n_2$, debris $n_3$ (i.e. all dead or apoptotic cells), chemoattractant $c_1$ (immune cells and SMCs attract to), non oxidized $c_2$ and oxidized LDL $c_3$. Let $\Omega\subset\mathbb{R}^d$, $d=2,3$ be the domain of the arterial wall, $\Gamma_1$ the boundary between the arterial wall and the lumen and $\Gamma_2$ the outer boundary of the arterial wall.

Let us suppose that for all $x\in\Omega$ and $t>0$ the following system holds:
\begin{eqnarray}\label{eq:01}
\partial_t n_1 &=& \nabla\cdot \left( \mu_1 \nabla n_1 -\chi(n_1,c_1,\chi_{11}^0,\chi_{11}^{th})\nabla c_1 - \chi(n_1,c_3,\chi_{13}^0,\chi_{13}^{th})\nabla c_3 \right) - d_1n_1,\\
\partial_t n_2 &=& \nabla\cdot \left( \mu_2 \nabla n_2 -\chi(n_2,c_1,\chi_{21}^0,\chi_{21}^{th})\nabla c_1 + \chi( n_2,n_1,\chi_{21}^0,\chi_{21}^{th})\nabla n_1 \right) - d_2n_2,\\
\partial_t n_3 &=& \nabla\cdot ( \mu_3 \nabla n_3 ) + d_1 n_1 + d_2 n_2 - F(n_3,c_3) n_1,\\
\partial_t c_1 &=& \nabla\cdot ( \nu_1 \nabla c_1 ) - \alpha_1 n_1 c_1 - \alpha_2 n_2 c_1 + f_1(n_3)n_3,\\
\partial_t c_2 &=& \nabla\cdot (\nu_2 \nabla c_2 ) - k c_2,\\
\partial_t c_3 &=& \nabla\cdot (\nu_3 \nabla c_3 ) + k c_2.
\end{eqnarray}
In our model we assume the motility coefficients $\mu_1$, $\mu_2$, $\mu_3$, $\nu_1$, $\nu_2$ and $\nu_3$ to be constant. The parameters $d_1$ and $d_2$ are also constant and describe the death rates of immune cells and SMCs. Chemoattractant is neutralized by immune cells and SMCs which is described by $\alpha_1$ and $\alpha_2$. The parameter $k$ describes how fast the native LDL becomes oxidized. 

The functions $\chi$ is called tactic sensitivity function.
We have chosen $\chi(x,y,a,b) = a \frac{x}{y+b}$ to mimic a high sensitivity of cells to the relative gradient $\frac{\nabla c}{c}$ of a chemoattractant (or other cells) $c$ on the one hand and a small penalization term to regularize the (chemo-)tactic movement for small concentrations $c$ on the other hand. A lot of other tactic sensitivity functions are possible as well. Our tactic sensitivity functions are defined by constants $\chi_{ij}^0$ and $\chi_{ij}^{th}$.

For a healthy immune system debris is degraded which is indicated by a general function $F>0$. We suppose $\gamma:=F<0$ to be constant indicating a diseased state. The function $f_1$ is a production term which is debris dependent. 
We allow LDL and immune cells to enter the arterial wall through the inner boundary and SMCs to enter through the outer arterial wall. The immune cell (SMC) inflow is triggered when a threshold $c_1^{*}$ ($c_1^{**}$) of chemoattractant is exceeded, i.e.
\begin{eqnarray}\label{eq:03}
	\partial_n n_1 &=& -\beta_1 H(c_1 -c_1^{*} ) \quad\quad \forall x\in \Gamma_1,\:t>0,\\
	\partial_n n_2 &=& -\beta_2 H(c_1 -c_1^{**} ) \quad\quad \forall x\in \Gamma_2,\:t>0,\\
	\partial_n c_2 &=& -\sigma \quad\quad\quad\quad\quad\quad\quad\quad \forall x\in \Gamma_{1,in},\:t>0,
\end{eqnarray}
with Heaviside function $H$ and a boundary $\Gamma_{1,in}\subset\Gamma_1$ for the inflow of LDL. Here, $\beta_1$, $\beta_2$ and $\sigma$ denote constant inflow rates for immune cells, SMCs and LDL, respectively. For all other boundary conditions we choose a no-flow condition.
The initial data is supposed to be given by $n_i(x,0)=n_i^0(x)$ and $c_i(x,0)=c_i^0(x)$, $i=1,2,3$, $x\in\Omega$.

Defining a vector $U := (n_1, n_2, n_3, c_1, c_2, c_3 )$ and functions $\mathcal{F}:\mathbb{R}^{6}\rightarrow\mathbb{R}^{6\times d}$, $\mathcal{A}:\mathbb{R}^{6}\rightarrow\mathbb{R}^{6\times 6}$ and $S:\mathbb{R}^{6}\rightarrow\mathbb{R}^{6}$ by
\begin{eqnarray*}
	\mathcal{F}(U)&:=&(\chi(n_1,c_1,\chi_{11}^0,\chi_{11}^{th})\nabla c_1 + \chi(n_1,c_3,\chi_{13}^0,\chi_{13}^{th})\nabla c_3,\\
      &&\chi(n_2,c_1,\chi_{21}^0,\chi_{21}^{th})\nabla c_1 - \chi( n_2,n_1,\chi_{21}^0,\chi_{21}^{th})\nabla n_1,0\ldots0),\\
	\mathcal{A}(U)&:=&\mathrm{diag}(\mu_1,\mu_2,\mu_3,\nu_1,\nu_2,\nu_3 ),\\
	S(U)&:=&-(d_1n_1,d_2n_2,-d_1 n_1 - d_2 n_2 + \gamma n_1,\alpha_1 n_1 c_1+\alpha_2 n_2 c_2 - f_1n_3, k c_2,-k c_2)
\end{eqnarray*}
equation \eqref{eq:01} can be written as 
\begin{equation}
	\partial_t U = - \nabla\cdot( \mathcal{F}(U) - \mathcal{A}(U)\nabla U) + S(U).
\end{equation}

\section{Discretization}
\label{seq:discretization}

The considered discretization is based on the Discontinuous Galerkin (DG) approach and 
implemented in \dunefem \cite{dedner2010generic}  a module of the
\dune framework \cite{bastian2008generic}.
The current state of development allows for simulation of convection dominated 
(cf.~\cite{limiter:11}) as well as viscous flow (cf.~\cite{brdar2012cdg2}).
We consider the CDG2 method from 
\cite{brdar2012cdg2} for various polynomial orders in space and 2nd (or 3rd) order 
in time for the numerical investigations carried out in this paper.

\subsection{Spatial Discretization}

The spatial discretization is derived in the following way. 
Given a tessellation $\grid$ of the domain $\Omega$ with 
$\cup_{K \in \grid} K = \Omega$ the 
discrete solution $\df$ is sought in the piecewise polynomial space 
\begin{equation}
\label{eqn:vspace}
    \phispace = \{\vecv\in L^2(\Omega,\mathbb{R}^{\nspec}) \; \colon
    \vecv|_{K}\in[\mathcal{P}_k(K)]^{\nspec}, \ K\in\grid\}
      \quad\textrm{for some}\;k \in \NNN, \nonumber
\end{equation}
where $\nspec$ is the number of species and $\mathcal{P}_k(K)$ is a space containing polynomials up to degree
$k$. 

\newcommand{\dual}[1]{\langle \basefct, #1 \rangle}%
We denote with $\Gamma_i$ the set of all intersections between two 
elements of the grid $\grid$ and accordingly with $\Gamma$ the set of all
intersections, also with the boundary of the domain $\Omega$. 
The following discrete form is not the most general but still
covers a wide range of well established DG methods. 
For all basis functions $\basefct \in \phispace$ we 
define 
\begin{equation}
\label{convDiscr}
\dual { \spcoper(\df) } := \dual{ \oper{K}_h(\df) } + \dual{ \oper{I}_h(\df) }
\end{equation}
with the element integrals  
\begin{eqnarray}
\label{eqn:elementint}
   \dual{ \oper{K}_h(\df) } &:=&
      \sum_{\elem \in \grid} \int_{\elem}
      \big( ( \mathcal{F}(\df) - \mathcal{A}(\vecU_h) \nabla \df ) : \nabla\basefct + \su
      \cdot \basefct \big),
\end{eqnarray}
and the surface integrals (by introducing appropriate numerical fluxes 
$\fluxF$, $\fluxA$ for the convection and diffusion terms, respectively) 
\begin{eqnarray}
\label{eqn:surfaceint}
   \dual{ \oper{I}_h(\df) } &:=&
      \sum_{e \in \Gamma_i} \int_e \big(
      \vaver{\mathcal{A}(\vecU_h)^T\nabla\basefct} : \vjump{\df} +
      \vaver{\mathcal{A}(\vecU_h)\nabla\df} : \vjump{\basefct} \big) \nonumber \\
    &-& \sum_{e \in \Gamma} \int_e \big( \fluxF(\df) - \fluxA(\df)\big) :
      \vjump{\basefct},      
\end{eqnarray}
where $\vaver{ \vecV } = \frac{1}{2}( \vecV^+ + \vecV^- )$ denotes the average and 
$\vjump{ \vecV } = (\nbold^+ \otimes \vecV^+  + \nbold^-\otimes \vecV^-) $ the jump of the
discontinuous function $\vecV\in V_h$ over element boundaries.
For matrices $\sigma,\tau\in\mathbb{R}^{m\times n}$ we use standard notation
$\sigma : \tau = \sum_{j=1}^m\sum_{l=1}^n\sigma_{jl}\tau_{jl}$. Additionally, for vectors
$\vecv \in \mathbb{R}^m,\vecw\in\mathbb{R}^n$, we define $\vecv\otimes\vecw\in\mathbb{R}^{m\times n}$
according to $(\vecv\otimes\vecw)_{jl}=\vecv_j \vecw_l$ for $1\leq j\leq m$, $1\leq l\leq n$.

The convective numerical flux $\fluxF$ can be any appropriate numerical flux known for
standard finite volume methods. 
For the results presented in this paper we choose $\fluxF$ to be the widely used
local Lax-Friedrichs numerical flux function. 

A wide range of diffusion fluxes $\fluxA$ can be found in the
literature, for a summary see \cite{arnold2002unified}.
We choose the CDG2 flux
\begin{eqnarray}
\fluxA(\vecV) := 2\liftfactor_e \big(\mathcal{A}(\vecV)\liftre(\vjump{\vecV})\big)|_{\Kminus}
\quad\mbox{for } \vecV\in V_h,
\end{eqnarray}
which was shown to be highly efficient for advection-diffusion equations (cf. \cite{brdar2012cdg2}). 
Based on stability results, we choose $\Kminus$ to be the element adjacent to the edge $e$ with the smaller
volume. $\liftre(\vjump{\vecV})\in [\phispace]^d$ is the lifting of the jump of $\vecV$ defined by
\begin{eqnarray}
    \int_\Omega \liftre(\vjump{\vecV}) : \taubold = -\int_e
    \vjump{\vecV} : \vaver{\taubold} \quad
    \mbox{for all}\;\taubold\in [\phispace]^d.
\end{eqnarray}
For the numerical experiments in this paper we use $\liftfactor_e= \frac{1}{2}\maxnuminterface$,
where $\maxnuminterface$ is the maximal number of intersections one element in the grid
can  have (cf. \cite{brdar2012cdg2}). We use triangular elements 
where $\liftfactor_e=1.5$ for all $\isec \in \Gamma$, and tetrahedral elements 
where $\liftfactor_e=2$ for all $\isec \in \Gamma$.

\subsection{Temporal discretization}
\label{TimeDisc}

The discrete solution $\df(t) \in \phispace$ 
has the form $\df(t,x) = \sum_i \vecU_i(t)\basefct_i(x)$.
We get a system of ODEs for the coefficients of $\vecU(t)$ which reads 
\begin{eqnarray}
  \label{eqn:ode}
  \vecU'(t) &=& f(\vecU(t),t)  \mbox{ in } (0,T]
\end{eqnarray}
with $f(\vecU(t),t) = M^{-1}\spcoper(\df(t),t)$, $M$ being the mass matrix which is in
our case block diagonal or even diagonal, depending on the choice of basis
functions. $\vecU(0)$ is given by the projection of $\vecU_0$ onto $\phispace$.

For the numerical results 
we have chosen Diagonally Implicit Runge-Kutta (DIRK) 
solvers of order $2$, $3$, or $4$ depending
on the polynomial order of the basis functions. The DIRK solvers are  
based on a Jacobian-free Newton-Krylov method (see \cite{knoll:04}).
The Krylov method is chosen to be GMRES. 
The implicit solver relies on a \textbf{matrix-free} implementation of the discrete operator 
$\spcoper$. In a follow-up paper we will compare this approach to a fully assembled
approach.

\section{Numerical Results}
In this section we present some benchmark tests for 2D and 3D focusing on parallelization and higher order DG schemes. Due to the lack of an exact solution $U$ we have computed the $L^2$-error between the discrete solution $U_h$ and a very fine, higher order solution $U_{h'}$. The quadrature order to compute $\|U_h-U_{h'}\|_{L^2(\Omega)}$ was chosen to be $2k+4$, where $k$ denotes the order of the scheme. All computations are done on an unstructured, tetrahedral mesh.

\subsection{A 2D numerical experiment with six species}

\begin{table}[t]
\small
\caption{Accuracy of the CDG2 scheme with 32 threads}
\begin{tabular}{ll|lll|lll|lll}
\hline\noalign{\smallskip}
 & & & linear & & & quadratic & & & cubic & \\
level & grid size & time$^a$ & $L^2$-error & EOC$^b$ & time$^a$ & $L^2$-error & EOC$^b$ & time$^a$ & $L^2$-error & EOC$^b$ \\
\hline
0 & 80    & 5.72E-1& 2.42E-3 & ---     & 2.00E0 & 2.18E-3 & ---     & 6.52E0 & 1.96E-3 & --- \\
1 & 320   & 5.56E0 & 2.10E-3 & 0.20311 & 2.33E1 & 1.82E-3 & 0.26650 & 8.63E1 & 1.50E-3 & 0.38074 \\
2 & 1280  & 3.98E2 & 1.82E-3 & 0.21263 & 2.09E2 & 1.34E-3 & 0.43315 & 8.22E2 & 9.26E-4 & 0.69823 \\
3 & 5120  & 3.33E3 & 1.39E-3 & 0.38944 & 2.21E3 & 7.92E-4 & 0.76429 & 9.12E3 & 4.32E-4 & 1.0993 \\
4 & 20480 & 3.01E4 & 8.28E-4 & 0.74208 & 2.10E4 & 2.94E-4 & 1.4284  & 8.02E4 & 8.77E-5 & 2.3024 \\
5 & 81920 & 2.67E5 & 3.21E-4 & 1.3659  & 1.93E5 & 7.26E-5 & 2.0193  & 6.96E5 & 2.33E-5 & 1.9122 \\
\hline
\end{tabular}
$^a$ total CPU time, $^b$ experimental order of convergence, 
\label{tab_CGD2}
\end{table}

\begin{figure}[t]
\includegraphics[scale=0.215]{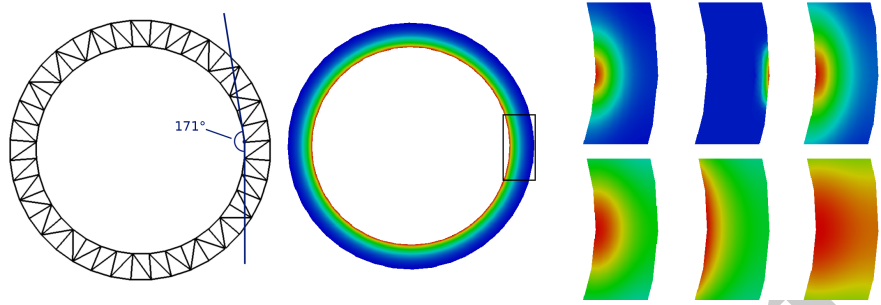}
\caption{Left: The coarsest grid for the EOC calculations containing 80 elements visualising a re-entrant corner (blue). The angle of $171^{\circ}$ stays fixed for all refinements. Middle: Initial distribution for the immune cells. Right: Solution for 6 species from left to right, up to down: Immune cells, SMCs, debris, chemoattractant, native LDL, oxidized LDL. (data visualisation: Paraview.)}
\label{fig:2d}
\end{figure}

\begin{figure}[t]
\includegraphics[scale=0.17]{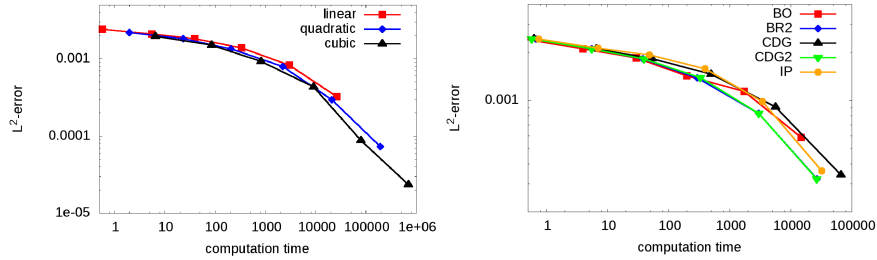}
\caption{Plot CPU time vs. $L^2$-error: left: 1st, 2nd and 3rd order CDG2 scheme, right: 1st order CDG, CDG2, Baumann-Oden (BO), Bassy-Rebay (BR2), interior penalty (IP) scheme (Visualisation of graphs: gnuplot.)}
\label{fig:eoc}
\end{figure}

$U_{h'}$ was calculated using the 4th order CDG2 scheme on a grid with 81,920 elements (refinement level $5$), i.e. 7,372,800 degrees of freedom. For each $h$-refinement of the grid we bisect the time step size. Results for linear, quadratic and cubic DG schemes can be seen in table \ref{tab_CGD2}. In figure \ref{fig:eoc} (left picture) we compare on a log-log scale the total CPU time of all threads with the $L^2$-error. Although the convergence rate is not as high as from the theory for parabolic problems, we see better rates for higher order schemes. We assume that re-entrant corners 
are responsible for the reduced convergence rates, see re-entrant corners in left picture of figure \ref{fig:2d}.

The right picture of figure \ref{fig:eoc} shows that the CDG2 is as good as the BR2 scheme and outperforms other DG schemes. 

\subsection{A 3D numerical experiment with three species}

For the 3D benchmark we simplify our model and do our simulation only for immune cells, debris and  chemoattractant. This reduces the considered model to 
\begin{eqnarray}\label{eq:02}
\partial_t n_1 &=& \nabla\cdot \left( \mu_1 \nabla n_1 -\chi(n_1,c_1,\chi_{11}^0,\chi_{11}^{th})\nabla c_1\right),\\
\partial_t n_3 &=& \nabla\cdot ( \mu_3 \nabla n_3 ) + d_1 n_1 + d_2 n_2 - F(n_3,c_3) n_1,\\
\partial_t c_1 &=& \nabla\cdot ( \nu_1 \nabla c_1 ) - \alpha_1 n_1 c_1 + f_1(n_3)n_3.
\end{eqnarray}
We cannot trigger the inflammation through an inflow of LDL anymore. Thus, we suppose that the inflammation is triggered by a local, high concentration of debris and keep all other boundary and initial data from the last section.

\begin{figure}[t]
\includegraphics[scale=0.16]{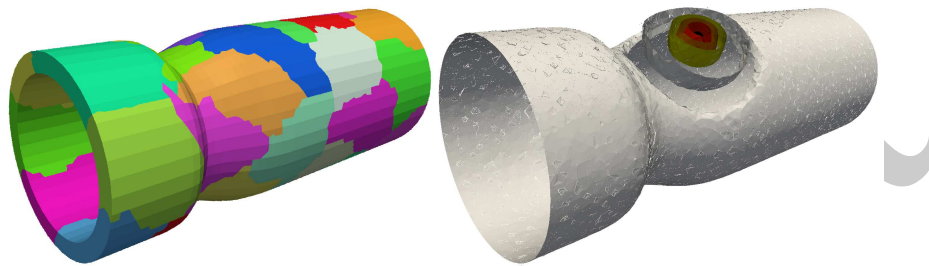}
\caption{3D cuff model. Left: Each colour denotes a processor in a parallel run with 32 processors, right: Isolines of the distribution of the chemoattractant after the inflammation has started}
\label{fig:3d}
\end{figure}

In the 3D benchmark we examine parallelization using MPI and present in table \ref{tab_parallel} strong scaling results for a third order CDG2 scheme on a grid with 113,549 elements and 13,625,880 degrees of freedom. Figure \ref{fig:3d} shows the distribution of the processors and a discrete solution of the chemoattractant calculated using first order CDG2.

\begin{table}[t]
\caption{CPU time for a parallel runs using the cubic CDG2 method for computation of
$10$ time steps.}
\begin{tabular}{p{2.5cm} p{1.0cm}p{1.0cm}p{1.0cm}p{1.0cm}p{1.0cm}p{1.0cm}}
\hline
processors      & 8      &   16   &   32   &   64   & 128    & 256 \\
CPU time in sec & $1177$ & $528$  & $277$  & $142$  & $75$   & $39$\\
speedup         &  ---   & $2.23$ & $4.29$ & $8.29$ & $15.7$ & $30.18$\\
\hline
\end{tabular}
\label{tab_parallel}
\end{table}

\section{Conclusion}

We have shown that Discontinuous Galerkin schemes are well suited for solving huge coupled reactive diffusion transport systems. Modern techniques, such as parallelization, help to handle large systems in an appropriate CPU time.
Furthermore, we have shown that it is possible to model the early stages of atherosclerotic plaque formation. A lot of more work needs to be done: In a future paper we will model the wall shear stress and some more species to understand later stages of atherosclerosis. 

\section*{Acknowledgement}
This work was supported by the Deutsche Forschungsgemeinschaft, Collaborative Research Center
SFB 656 ``Cardiovascular Molecular Imaging'', project B07, M{\"u}nster, Germany. The scaling results were produced using the super computer Yellowstone
(ark:/85065/d7wd3xhc) provided by NCAR's Computational and Information Systems
Laboratory, sponsored by the National Science Foundation.
Robert Kl\"ofkorn is partially funded by 
the DEO program BER under award DE-SC0006959.

\bibliographystyle{unsrt}
\bibliography{girke2014efficient}

\end{document}